\documentclass{article}
\usepackage[ansi]{umlaute}
\usepackage{color}
\usepackage{colordvi}
\usepackage{amsmath, amssymb, amsthm}
\usepackage{graphics,eepic}

\newtheorem{lemma}{Lemma}
\newtheorem{proposition}{Proposition}
\newtheorem{theorem}{Theorem}
\theoremstyle{definition}

\theoremstyle{remark}

\newtheorem{remark}{Remark}

\newcounter{marginnote}
\definecolor{addedcolor}{rgb}{0,0,0.4}
\definecolor{removedcolor}{rgb}{0.4,0,0}
\definecolor{marginnumbercolor}{rgb}{0.6,0.6,0}

\newcommand{\inserted}[1]{{\color{addedcolor}#1}}

\newcommand{\poly}{\ensuremath{\mathcal P}}
\newcommand{\binomial}[2]{\ensuremath{\left( \begin{matrix}#1 \\ #2 \end{matrix} \right)}}
\newcommand{\inner}[2]{\ensuremath{\left \langle #1 , #2 \right \rangle}}
\newcommand{\innerP}[2]{\ensuremath{\left \langle #1 , #2 \right \rangle_{\poly}}}
\newcommand{\dd}{\ensuremath{\mathrm{d}}}
\newcommand{\diag}{\ensuremath{\mathrm{diag}}}

\newcommand{\res}{\ensuremath{\mathrm{Res}}}

\newcommand{\sym}{\ensuremath{\mathcal{S}_n}}
\newcommand{\Span}{\ensuremath{\mathrm{Span}}}
\newcommand{\vol}{\ensuremath{\mathrm{Vol}}}

\newcommand{\real}{\ensuremath{\mathrm{Re}}}
\newcommand{\imag}{\ensuremath{\mathrm{Im}}}

\newcommand{\R}{\mathbb{R}}
\newcommand{\N}{\mathbb{N}}
\newcommand{\On}{\mathbb{O}_n}
\newcommand{\ev}{\mbox{eval}_1}

\author{Jean-Pierre Dedieu\thanks{
Institut de Mathématiques.
Université Paul Sabatier. 
31062 Toulouse Cedex 9.
France.
e-mail : jean-pierre.dedieu@math.ups-tlse.fr
} \and
Gregorio Malajovich\thanks{Departamento de Matem\'atica Aplicada,
Universidade Federal do Rio de Janeiro. Caixa Postal 68530, CEP 21945-970,
Rio de Janeiro, RJ, Brasil. e-mail: gregorio@ufrj.br}}
\title{On the number of minima of a random polynomial\thanks{
This work was sponsored by the International Cooperation Agreement 
Brazil-France. J-P Dedieu was supported by the ANR Gecko.
G. Malajovich was also supported by Brazilian CNPq
grants 304504/2004-1 and 472486/2004-7, and by FAPERJ 
(Fundação Carlos Chagas de Amparo
à Pesquisa do Estado do Rio de Janeiro) grant 
E26/170.734/2004/Edital Primeiros Projetos.
}
}
\date{February 9, 2007}

\begin{document}
\maketitle
\begin{abstract}
We give an upper bound in $O(d^{(n+1)/2})$ for the number
of critical points of a normal random polynomial with degree at most $d$ and $n$ variables. 
Using the large deviation principle for the spectral value of large random matrices we obtain the bound 
$$O\left(\exp(-\beta n^2 + \frac{n}{2}\log(d-1))\right)$$
($\beta$ is a positive constant independent on $n$ and $d$) 
for the number of minima of such a polynomial. 
This proves that most normal random polynomials of fixed degree have only saddle points.
Finally, we give a closed form expression for the
number of maxima (resp. minima) of a random univariate polynomial,
in terms of hypergeometric functions.
\end{abstract}

\section{Introduction}
We consider a random polynomial $f$ over the reals with $n \ge 1$ variables and
degree at most $d \ge 2$. The problem is to compute, on the average, its number of
critical points (the number of real roots of the system $Df(x) = 0$), and its number
of local minima. Since a generic polynomial has only nondegenerate stationnary
points, this last number is also given by the real roots of the system $Df(x) = 0$
such that $D^2f(x)$ is positive definite. This {reduces} our problem to the
computation of the number of real roots of a polynomial system under certain
constraints. 

Generally speaking, let $F=(F_1, \ldots , F_n)$ be a random system of real
polynomial equations with $n$ variables and degree $F_i \le d_i.$ Let $N^F(U)$
denotes the number of zeros of the system $F(x)=0$ lying in the subset $U \subset
\R^n$ and $N^F(\R^n ) = N^F$. Little is known on the distribution of the random
variable $N^F(U)$. A classical result in the case of one polynomial of one variable
is given by Kac \cite{kac1}, \cite{kac2}, who gives the asymptotic value
$$E(N^F) \approx \frac{2}{\pi} \log d $$
as $d$ tends to infinity when the coefficients of $F$ are Gaussian centered
independent random variables with variances equal to $1$. But, when the variance of
the $i-$th coefficient is equal to $d \choose i$ (Weyl's distribution), we have (see
Bogomolny-Bohias-Leboeuf \cite{bog} and also Edelman-Kostlan \cite{ede})
$$E(N^F) = \sqrt{d}.$$

In 1992, Shub and Smale extended this result to a real polynomial system $F$ where 
$$F_i(x_1, \ldots , x_n) = \sum_{\alpha_1 + \ldots + \alpha_n \le
d_i}a_{i,\alpha}x_1^{\alpha_1} \ldots x_n^{\alpha_n},$$
when the coefficients $a_{i,\alpha}$ are Gaussian centered independent random
variables with variances equal to 
$$\binomial{d}{\alpha} = \frac{d_i !}{\alpha_1! \ldots \alpha_n! (d_i - \alpha_1
\ldots - \alpha_n)!}$$
(see Kostlan \cite{kos} on this distribution and its properties). Their result is
$$E(N^F) = \sqrt{d_1 \ldots d_n}$$
that is the square root of the B\'ezout number of the system. 

A general formula for the expected value of $N^F(U)$ when the random functions
$F_i$, $1 \le i \le n,$ are stochastly independent and their law is centered and
invariant under the isometries of $\R^n$ can be found in Aza\"is-Wschebor
\cite{aza}. This includes the Shub-Smale formula as a special case. 

This result has also been extended by Rojas \cite{roj} to {multi-homogeneous}
polynomial systems, and 
then partially by Malajovich and Rojas~\cite{mal} to sparse 
polynomial systems. 

Wschebor in \cite{wsc} studies the moments of $N^F$ and Armentano-Wschebor
\cite{arm} consider random systems of equations of the type $P_i(x) + X_i(x)$, $1
\le i \le n,$ $x \in \R^n,$ where the $P_i's$ are non-random polynomials (the
signal) and the $X_i's$ are independent Gaussian random variables (the noise). 

Notice a major difference between these studies and the case considered here: the
$n$ equations of the system $Df(x)=0$ are not independent!

Through this paper we denote by $\poly = \poly_{d,n}$ the space of degree at most $d$, $n$-variate polynomials with
real coefficients. This space is endowed with the inner product:
\[
\innerP{f}{g} = \sum_{|\alpha| \le d} \binomial{d}{\alpha}^{-1} f_\alpha g_\alpha 
\]
where $\alpha = (\alpha_1, \ldots , \alpha_n) \in \N^n$ is a multi-integer, $|\alpha| = \alpha_1 + \ldots + \alpha_n$, 
$$f(x) = \sum_{ |\alpha| \le d}f_{\alpha} x_1^{\alpha_1} \ldots x_n^{\alpha_n} = \sum_{ |\alpha| \le d}f_{\alpha} x^{\alpha},$$
and again
$$\binomial{d}{\alpha} = \frac{d !}{\alpha_1! \ldots \alpha_n! (d - |\alpha|)!}.$$
We make $\poly$ a probability space in considering the probability measure
\[
\frac{1}{\sqrt{2 \pi}^{\dim \poly}} e^{-\|f\|_{\poly}^2 / 2} \dd \poly
= \frac{1}{\sqrt{2 \pi}^{\dim \poly}} e^{-\|f\|_{\poly}^2 / 2}
\bigwedge_{|\alpha| \le d} \binomial{d}{\alpha}^{-1/2}{df_\alpha}
\]
i.e. a random polynomial has here Gaussian centered independent random coefficients with variances equal to $\binomial{d}{\alpha}$. 

Let $\mathcal S_n$ be the space of $n \times n$ real symmetric matrices, endowed with the Frobenius inner product $\langle R, S \rangle = \mbox{Trace}(R^T S)$ and its induced
norm $$\|S\|^2 = \sum_{1 \le i,j \le n} S_{ij}^2.$$ The {\em Gaussian Orthogonal Ensemble} is the space $\mathcal S_n$ together with the
probability measure
\[
\frac{e^{-\|S\|^2 / 2}}{(2 \pi)^{n(n+1)/4}}dS = \frac{e^{-\|S\|^2 / 2}}{2^{n/2} \pi^{n(n+1)/4}}\bigwedge_{1 \le i \le j \le n}{dS_{ij}}.
\]
Thus, the entries of a matrix in $\mathcal S_n$ are independent Gaussian random variables with mean $0$ and variance $1$ for a diagonal entry, and mean $0$ and variance $1/2$ for a non-diagonal entry.

Our first main result is the following:

\begin{theorem} \label{general}
Let $C_{d,n}$ denote the expected number of critical points of a random polynomial of degree at most $d$ in $n$ variables, and $E_{d,n}$ the expected number of minima. Let $P_n$ be the probability that a matrix in the Gaussian Orthogonal Ensemble is positive definite.
Then, for every $n \ge 2$,
$$C_{2,n} = 1 \mbox{ and }  E_{2,n} = P_n,$$
and for $d \ge 3$ 
\[ C_{d,n} \le \sqrt{\frac{2}{d}} (d-1)^{(n+2)/2} \mbox{ and } E_{d,n} \le \sqrt{\frac{2}{d}} (d-1)^{(n+2)/2} P_n.
\]
When $n=1$ one has
\[
C_{d,1} = 2E_{d,1} = \frac{2\sqrt{d-1}}{\pi} \int_0^\infty \frac{\sqrt{d(d-1)r^4 + 2dr^2 + 2}}{(dr^2 + 1)(r^2 + 1)} dr \le 1 + \sqrt{d-2}.
\]
Moreover, when $d \rightarrow \infty$, 
$$\frac{C_{d,1}}{1 + \sqrt{d-2}} \rightarrow 1.$$
\end{theorem}

Let $P_n$ be the probability that a matrix in the Gaussian Orthogonal Ensemble $GOE(n)$ is positive definite:
$$P_n = \int_{\sym^{++}} \frac{e^{-\|S\|^2 / 2}}{2^{n/2} \pi^{n(n+1)/4}}\bigwedge_{1 \le i \le j \le n}{dS_{ij}}.$$
Via the change of variable $S = Q \Lambda Q^T$ with $Q \in \On$ and $\Lambda = \diag (\lambda_1 \ge \ldots \ge \lambda_n \ge 0)$ one has
$$P_n = \frac{\mbox{Vol }\On}{2^n} \int_{\R^n_>}\prod_{i<j}(\lambda_i - \lambda_j) \frac{e^{-\| \lambda \|^2 / 2}}{(2 \pi )^{n(n+1)/4}} d\lambda $$
where $\lambda \in \R^n_>$ if and only if $\lambda_1 > \ldots > \lambda_n > 0$ and 
$$\mbox{Vol }\On = \frac{2^{n(n+3)/4}\Gamma(1/2)^{n(n+1)/2}}{\prod_{j=1}^n \Gamma((n-j+1)/2)}$$
(see Mehta \cite{meh} for the description of $P_n$ as an integral over $\R^n$ and Federer \cite{fed} for the volume of the orthogonal group). The following values are easy to obtain
$$P_1 = \frac{1}{2},\ P_2 = \frac{2-\sqrt{2}}{4},\ P_3 = \frac{\pi - 2\sqrt{2}}{4 \pi}.$$ 
$P_3$ was computed by Carlos Beltr\'an.

Using the large deviation principle for the spectral value of large random matrices (see the appendix at the end of this paper) we see that the asymptotic value of $P_n$ for large values of $n$ satisfies
$$\limsup_{n \rightarrow \infty} \frac{1}{n^2}\log(P_n) \le - \alpha$$
where $\alpha $ is a positive constant independent on $n$. Thus, there exist two positive constants $\beta$ and $\gamma$ such that, for every $n \ge 1$, 
$$P_n \le \gamma e^{- \beta n^2}.$$
This gives our second main theorem:

\begin{theorem} There exist two positive constants $\beta$ and $K$ such that for every $n$ and $d$ the number of minima of a random polynomial satisfies
$$E_{d,n} \le K e^{- \beta n^2 + \frac{n}{2}\log(d-1)}.$$
\end{theorem}

\begin{remark} This is a quite surprising result : it shows that most of random polynomials of reasonable degree have only saddle points.
\end{remark}

We want to thank here Alice Guionnet, Manjunath Krishnapur and Balint Virag who introduced us in the world of large deviation delicacies. 

\section{The space of $n$-variate polynomials}

The inner product space $\poly, \innerP{\cdot}{\cdot}$ has several interesting properties resumed in the following 

\begin{lemma} \label{lem-1} \begin{enumerate} \item It admits the {\em reproducing kernel}
$K(z,x) = \left(1 + \inner{z}{x}\right)^d$:
\begin{equation}\label{rep1}
f(x) = \innerP{K(.,x)}{ f }
\end{equation}
for any $x \in \R^n$ and $f \in \poly$.

\item It has a representation {formula} for the derivatives: for any integer $k \ge 1$ and $x,u_1, \cdots ,u_k \in \R^n$ we have
\begin{equation}\label{rep2}
D^kf(x)(u_1,\cdots ,u_k)  = \innerP{K_k(.,x,u_1,\cdots ,u_k)}{ f },
\end{equation}
with 
\begin{equation}\label{rep3} K_k(z,x,u_1,\cdots ,u_k) = D_x^kK(z,x)(u_1,\cdots ,u_k) = \end{equation}
$$d \cdots (d-k+1)\inner{z}{u_1} \cdots \inner{z}{u_k}\left(1 + \inner{z}{x}\right)^{d-k}.$$

\item This scalar product is orthogonally invariant:
\begin{equation}\label{rep4} \innerP{f \circ U}{g \circ U} = \innerP{f}{g}\end{equation}
for any $f,g \in \poly$ and the orthogonal transformation $U \in \On$.
\end{enumerate}
\end{lemma}

\begin{proof} The first two formulas are well known and easily obtained via a direct computation. For the orthogonal invariance
see \cite{blu}, section 12.1, or \cite{kos}. 
\end{proof}

A second interest of Weyl's distribution for polynomials is due to the following identity: let $f(x) = x^T S x$ (here $S$ is a symmetric $n \times n$ matrix) be a homogeneous degree $2$ polynomial, then $\| f \|_\poly = \| S \|$. This is the reason why 

\begin{proposition} \label{pro-0} $C_{2,n}= 1$ and $E_{2,n} = P_n$. \end{proposition}

\proof Since a generic degree $2$ polynomial has only one critical point we have $C_{2,n}= 1$. Given $f \in \poly_{2,n}$ we can write it 
$$f(x) = \alpha + \sum_{1\le i \le n} b_i x_i + \sum_{1\le i \le n} a_{ii} x_i^2 + \sum_{1\le i < j \le n} a_{ij} x_i x_j.$$
One has
$$\| f \|_\poly^2 = \alpha^2 + \frac{1}{2}\sum_{1\le i \le n} b_i^2 + \sum_{1\le i \le n} a_{ii}^2 + \frac{1}{2} \sum_{1\le i < j \le n} a_{ij}^2$$
so that
$$E_{2,n} = \int_{D^2f(0) > 0} \frac{e^{-\| f \|_\poly^2/2}}{2^{n(n+1)/4} (2\pi)^{(n+1)(n+2)/4}}d\alpha db da = $$
$$\int_{D^2f(0) > 0} \frac{e^{-\left(\sum_{i} a_{ii}^2 + \frac{1}{2} \sum_{i < j} a_{ij}^2\right)/2}}{2^{n(n-1)/4} (2\pi)^{n(n+1)/4}}da.$$
To compute this last integral we let $S = \frac{1}{2}D^2f(0)$; this gives
$$E_{2,n} = \int_{S > 0} \frac{e^{-\left\|S\right\|^2/2}}{2^{n(n-1)/4} (2\pi)^{n(n+1)/4}}2^{n(n-1)/2}dS = P_n.$$
\qed

\section{An integral formulation}

Let us define 
$$\ev : \poly \times \R^n \rightarrow \R^n, \ \ev(f,x) = Df(x).$$
The incidence variety for real critical points of a polynomial is
defined by
\[
V = \left\{ (f,x) \in \poly \times \mathbb R^n \ : \ \ev(f,x) = 0 \right\}
\ .
\]
The derivative of $\ev$ is given by
$$D \ev (f,x) (\dot f, \dot x) = D\dot f(x) + D^2f(x)\dot x$$
for any $f,\dot f \in \poly$ and $x, \dot x \in \R^n$. Since this derivative is onto, $V$ is a submanifold and its dimension is
$$\dim V = \dim \poly = {n+d \choose d}.$$
The tangent space at $(f,x) \in V$ is given by
$$T_{(f,x)}V = \ker D \ev (f,x) = \left\{(\dot f, \dot x) \in \poly \times \mathbb R^n \ : \ D\dot f(x) + D^2f(x)\dot x = 0 \right\}.$$

The restriction $\pi_2 : V \rightarrow \R^n$ of the projection $\poly \times \mathbb R^n \rightarrow \mathbb R^n$ is surjective and is also a regular map because for any $(f,x) \in V$ the derivative $D\pi_2 (f,x) : T_{(f,x)}V \rightarrow \R^n$ is surjective. The fiber of $\pi_2$ above $x \in \R^n$ 
$$V_x  = \left\{ (f,x) \in \poly \times \mathbb R^n \ : \ \ev(f,x) = 0 \right\}$$
is isomorphic to a $\dim \poly - n$ linear space. $V_x$ is equipped with the volume form inherited from the induced metric. 

The restriction $\pi_1 : V \rightarrow \poly$ of the projection $\poly \times \mathbb R^n \rightarrow \poly$ is a smooth map. A given $f \in \poly$ is a regular value of $\pi_1$ when either $f$ has no critical point or when, for any $x$ such that $(f,x) \in V$,  $D\pi_1 (f,x) : T_{(f,x)}V \rightarrow \poly$ is surjective. This last condition is satisfied when the second derivative $D^2f(x)$ is an isomorphism which is the generic situation: 
$$\Sigma '  = \left\{ (f,x) \in V \ : \ \det D^2f(x) = 0 \right\}$$
is a submanifold in $V$ and $\dim \Sigma ' < \dim V$. Thus $\Sigma '$ and its image $\Sigma = \pi_1 (\Sigma ')$ have zero measure and we may ignore them. For any $(f,x) \in V \setminus \Sigma '$ and any $\dot f \in \poly$ we have $D\pi_1(f,x)(\dot f , \dot x ) = \dot f$ for $\dot x = - D^2f(x)^{-1}D\dot f (x)$ and the fiber above $f$ 
$$V_f  = \left\{ (f,x) \in \poly \times \mathbb R^n \ : \ \ev(f,x) = 0 \right\}$$
consists in a finite number of points.

Given $(f,x) \in V \setminus \Sigma '$ we are in the context of the implicit function theorem that is $V$ is locally around $(f,x)$ the graph of the function
$$G = \pi_2 \circ \pi_1^{-1}$$
where $\pi_1^{-1}$ is the local inverse of $\pi_1$ such that $\pi_1^{-1}(f) = (f,x)$. Since the graph of $DG(f)$ is the tangent space $T_{(f,x)}V$ we get
\begin{equation}\label{rep5} DG(f) \dot f =  - D^2f(x)^{-1}D\dot f (x) \end{equation}
for any $\dot f \in \poly$. 

Like in \cite{blu} section 13.2, theorem 3, we have the following

\begin{proposition} \label{pro-1}
Let $U$ be a measurable subset of $V$. Let us denote by $\# (f,U)$ the number of pairs $(f,x) \in U$ and by $E_U$ the expectation of 
$\# (f,U)$ when $f$ is taken at random:
\begin{equation}\label{rep6} 
E_U = 
\int_{\poly}  \# (f, U)
\ \frac{e^{-\|f\|_{\poly}^2 / 2}}{(2 \pi)^{\dim \poly/2}}  \dd \poly .
\end{equation}
With these notations, one has 
\begin{equation}\label{rep7} 
E_U
\int_{\mathbb R^{n}}
\dd x \ 
\int_{V_x \cap U}
\det (DG(f)DG(f)^*)^{-1/2}  
\ \frac{e^{-\|f\|_{\poly}^2 / 2}}{(2 \pi)^{\dim \poly/2}} \dd V_x.
\end{equation}
\end{proposition}

\begin{remark} In our context two sets are of particular interest: $U = V$ to compute the average number of critical points of a polynomial $C_{d,n}$, and $U=V_+$ with
$$V_+ = \left\{(f,x) \in  \poly \times \mathbb R^n \ : \ Df(x) = 0 \mbox{ and } D^2f(x) > 0 \right\}$$ 
(here $>0$ means {\it positive definite}) for the average number of local minima $E_{d,n}$.
\end{remark}

We have now to compute the determinant appearing in equation \ref{rep7}. This is done in the following

\begin{proposition} \label{pro-2} Under the notations above
\begin{equation}\label{rep8} \det (DG(f)DG(f)^*) = d^n(1 + \left\| x \right\|^2)^{n(d-1)-1} (1 + d \left\| x \right\|^2) \left| \det D^2f(x) \right|^{-2}.
\end{equation}
\end{proposition}

\begin{proof} Let us denote $D\dot f (x)  =  D_x \dot f $. Since $DG(f) \dot f =  - D^2f(x)^{-1}D_x \dot f $ and since $D^2f(x)$ is symmetric,
we get 
$$DG(f)DG(f)^* = D^2f(x)^{-1} D_x D_x^* D^2f(x)^{-1}$$
so that
\begin{equation}\label{rep9} \det (DG(f)DG(f)^*) = \det (D_x D_x^*) \left| \det D^2f(x) \right|^{-2}. \end{equation}
To compute $\det (D_x D_x^*)$ we use the representation formula for the derivative (equation \ref{rep2}) with $k=1$. Let us denote by $e_i$, $1 \leq i \leq n,$ the canonical basis in $\R^n$. Then, for any $\dot f \in \poly$, 
$$D_x\dot f = \sum_i e_i \left\langle K_1(.,x,e_i), \dot f \right\rangle_\poly$$
so that, with $\dot x \in \R^n$, $\dot x = \sum_i \dot x_i e_i$, 
$$\left\langle D_x^* \dot x, \dot f \right\rangle_\poly = \left\langle \dot x, D_x \dot f \right\rangle \left\langle \dot x, \sum_i e_i \left\langle K_1(.,x,e_i), \dot f \right\rangle_\poly \right\rangle \sum_i \dot x_i \left\langle K_1(.,x,e_i), \dot f \right\rangle_\poly.$$
Thus, we get 
$$D_x^* \dot x = \sum_i \dot x_i K_1(.,x,e_i)$$
and consequently
$$D_x D_x^* \dot x = \sum_i e_i \left\langle K_1(.,x,e_i), \sum_j \dot x_j K_1(.,x,e_j) \right\rangle_\poly =$$
$$\left. \sum_i e_i \frac{\partial}{\partial z_i}\left(\sum_j \dot x_j d\left\langle z,e_j \right\rangle
\left( 1 + \left\langle z,x \right\rangle\right)^{d-1}
\right)\right|_{z=x} =$$
$$\sum_{i,j} e_i \dot x_j \times
\left\{
\begin{array}[pos]{ll}
	d(d-1)x_ix_j (1 + \left\| x \right\|^2)^{d-2}& \mbox{ if } i \ne j\\
	d(d-1)x_i^2 (1 + \left\| x \right\|^2)^{d-2} + d(1 + \left\| x \right\|^2)^{d-1}& \mbox{ if } i = j\\
\end{array}
\right. $$
which correspond to the matrix 
$$d(d-1)(1 + \left\| x \right\|^2)^{d-2}xx^T + d(1 + \left\| x \right\|^2)^{d-1}I_n.$$
Its eigenvectors are $x$ and any nonzero vector in the orthogonal subspace $x^\perp$. The corresponding eigenvalues are
$$d(d-1)(1 + \left\| x \right\|^2)^{d-2}\left\| x \right\|^2 + d(1 + \left\| x \right\|^2)^{d-1} = d(1 + \left\| x \right\|^2)^{d-2}(1 + d \left\| x \right\|^2)$$
with multiplicity $1$, and
$$d(1 + \left\| x \right\|^2)^{d-1}$$
with multiplicity $n-1$ so that
$$\det D_x D_x^* = d^n(1 + \left\| x \right\|^2)^{n(d-1)-1} (1 + d \left\| x \right\|^2).$$
Our proposition combines this value and equation \ref{rep9}. 
\end{proof}

If we combine propositions \ref{pro-1} and \ref{pro-2} we obtain the following integral formulation

\begin{proposition} \label{pro-3}
Let $U$ be a measurable subset of $V$. One has
\begin{equation}\label{rep61}
E_U =\int_{\R^n}dx \int_{V_x \cap U} \frac{\left| \det D^2f(x) \right|}{d^{n/2} (1+\left\|x\right\|^2)^{(n(d-1)-1)/2} (1+d\left\|x\right\|^2)^{1/2}}
\ \frac{ e^{-\|f\|_{\poly}^2/2}}{(2 \pi)^{\dim \poly/2}}  \dd V_x .
\end{equation}
\end{proposition}

An action of the orthogonal group $\On$ on $\poly \times \R^n$ is defined by 
$$(Q,f,x) \in \On \times \poly \times \R^n \rightarrow (f \circ Q, Q^T x) \in \poly \times \R^n.$$
This action leaves the incidence variety $V$ invariant and also the scalar product $\left\langle .,. \right\rangle_\poly$ (lemma \ref{lem-1}). For this reason, when the measurable set $U$ is itself invariant, the integral on $V_x \cap U$ in proposition \ref{pro-3} only depends on $r = \left\| x \right\|$. Thus, taking spherical coordinates in $\R^n$, we get:

\begin{proposition} \label{pro-4}
Let $U$ be a measurable subset of $V$ invariant under the action of $\On$ (for any $(Q,f,x) \in \On \times U$ we have $(f \circ Q, Q^T x) \in U$). Under this condition
{
\[
E_U
\frac{\alpha_n}{d^{n/2}}
\int_{0}^{\infty}
\frac{r^{n-1}\dd r}{
R^{(d-1)n-1} (dr^2+1)^{1/2}
}\int_{V_{r e_1} \cap U}| \det D^2 f (re_1) |
\frac{  e^{-\|f\|_{\poly}^2 / 2} }
{(2 \pi)^{\dim V_{r e_1}/2}}\dd V_{r e_1}
\]
where $\alpha_n = \frac{ \vol S^{n-1}}{(2 \pi)^{n/2}}
=\frac{2}{2^{n/2} \Gamma(n/2)}$,  $R=\sqrt{r^2 + 1}$ and $re_1^T = (r,0, \ldots , 0)$.
}

\end{proposition}

\begin{remark} The measurable sets considered here: $U = V$ and $U = V_+ = \left\{(f,x) \in V \ : \ D^2f(x) > 0 \right\}$, are clearly invariant under the action of $\On$.
\end{remark}
\section{The inner integral}

Our objective is now to compute the integral over $V_{r e_1} \cap U$ appearing in proposition $\ref{pro-4}$. 

Let $D^2: V_{r e_1} \rightarrow \sym$ denote the operator $f \mapsto D^2 f(re_1)$.
We would like to compute its pseudo-inverse $\Psi: \sym \rightarrow
(\ker D^2)^{\perp}$. This means that $\Psi$ is the minimum norm right inverse 
of $D^2$ ($D^2 \circ \Psi = \mbox{id}_{\sym}$).

This will allow us to `integrate out' $\ker D^2$:
\begin{equation}\label{intout}
\int_{V_{r e_1} \cap U} | \det D^2 f |
\frac{e^{-\| f \|_\poly^2 / 2}}{(2 \pi)^{\dim V_{r e_1}/2}}  \dd V_{r e_1}
= \end{equation}
$$
\int_{D^2(\inserted{V_{r e_1} \cap }U)}\left| \det S \right| \left| \det \Psi^*\Psi \right|^{1/2}  
\frac{e^{-\|\Psi(S)\|_\poly^2/2}}{(2 \pi)^{\dim \sym / 2}} 
\  \dd S .
$$
To compute $\Psi(S)$ and $\left| \det \Psi^*\Psi \right|$ we need the following lemma:

\begin{lemma} \label{lem-2} Let us denote  
\begin{itemize} 
\item $e_i$, $1 \le i \le n$, the canonical basis in $\R^n$,
\item $\partial_{e_i} = K_1(z,re_1,e_i)$,
\item $\partial_{e_ie_j} = K_2(z,re_1,e_i,e_j)$,
\item $R=\sqrt{1+r^2}$. 
\end{itemize}
Then,
\begin{enumerate}
\item 
$\innerP{\partial_{e_1}}{\partial_{e_1}}
d(1+dr^2)R^{2d-4} $
\item If $i \ne 1$, then $\innerP{\partial_{e_i}}{\partial_{e_i}}
= d R^{2d-2}$
\item  If $i \ne j$, then $\innerP{\partial_{e_i}}{\partial_{e_j}} = 0$
\item  $\innerP{\partial_{e_1}}{\partial_{e_1e_1}} = 
d(d-1)(dr^2+ 2)rR^{2d-6}$
\item  If $(i,j,k) \ne (1,1,1)$, then $\innerP{\partial_{e_j}}{\partial_{e_ie_k}} = 0$
\item $\innerP{\partial_{e_1e_1}}{\partial_{e_1e_1}} = d(d-1)\left( d(d-1)r^4 + 4(d-1)r^2+2 \right)R^{2d-8}$
\item If $k \ne 1$, then $\innerP{\partial_{e_1e_k}}{\partial_{e_1e_k}} = 
d(d-1) ((d-1)r^2 + 1) R^{2d-6}$
\item If $i \ne 1$ and $k \ne 1$ , then $\innerP{\partial_{e_ie_k}}{\partial_{e_ie_k}} =$ {$(1+\delta_{ik})$} $d(d-1) R^{2d-4}$ ($\delta_{ik}$ is the Kronecker symbol),
\item If $\left\{ i,k \right\} \ne \left\{ j,l \right\}$, then $\innerP{\partial_{e_ie_k}}{\partial_{e_je_l}} = 0$
\end{enumerate}
\end{lemma}

\begin{proof} It is a consequence of the representation formulas given in lemma \ref{lem-1}:
\begin{itemize}\item $\left\langle  \partial_{e_1}, \partial_{e_1}\right\rangle_\poly = \left\langle  K_1(.,re_1,e_1), K_1(.,re_1,e_1)\right\rangle_\poly = \frac{\partial}{\partial z_1} K_1(z,re_1,e_1)\left|_{z=re_1}\right. = \frac{\partial}{\partial z_1} dz_1(1 + rz_1)^{d-1}\left|_{z=re_1}\right. = d(1+r^2)^{d-2}(1+dr^2),$
\end{itemize}
and similarly
\begin{itemize}\item $\left\langle  \partial_{e_i}, \partial_{e_i}\right\rangle_\poly = \frac{\partial}{\partial z_i} K_1(z,re_1,e_i)\left|_{z=re_1}\right. = \frac{\partial}{\partial z_i} dz_i(1 + rz_1)^{d-1}\left|_{z=re_1}\right. = d(1+r^2)^{d-1},$
\item $\left\langle  \partial_{e_i}, \partial_{e_j}\right\rangle_\poly = \frac{\partial}{\partial z_i} K_1(z,re_1,e_j)\left|_{z=re_1}\right. = \frac{\partial}{\partial z_i} dz_j(1 + rz_1)^{d-1}\left|_{z=re_1}\right. = 0$ when $i \ne j$,
\item $\left\langle  \partial_{e_1}, \partial_{e_1e_1}\right\rangle_\poly = \frac{\partial}{\partial z_1} K_2(z,re_1,e_1,e_1)\left|_{z=re_1}\right. = \frac{\partial}{\partial z_1} d(d-1)z_1^2 (1 + rz_1)^{d-2}\left|_{z=re_1}\right. = 
 d(d-1)r(2+dr^2)(1+r^2)^{d-3},$
\item $\left\langle  \partial_{e_j}, \partial_{e_ie_k}\right\rangle_\poly = \frac{\partial}{\partial z_j} d(d-1)z_iz_k (1 + rz_1)^{d-2}\left|_{z=re_1}\right. = 0$ when $(i,j,k) \ne (1,1,1),$
\item $\left\langle  \partial_{e_1e_1}, \partial_{e_1e_1}\right\rangle_\poly = \frac{\partial^2}{\partial z_1^2} d(d-1)z_1^2 (1 + rz_1)^{d-2}\left|_{z=re_1}\right. = d(d-1) (1+r^2)^{d-4} (2 + 4(d-1)r^2 + d(d-1)r^4),$
\item $\left\langle  \partial_{e_1e_k}, \partial_{e_1e_k}\right\rangle_\poly = \frac{\partial^2}{\partial z_1z_k} d(d-1)z_1z_k (1 + rz_1)^{d-2}\left|_{z=re_1}\right. = d(d-1) (1+r^2)^{d-3} (1 + (d-1) r^2),$
\item $\left\langle  \partial_{e_ie_k}, \partial_{e_ie_k}\right\rangle_\poly = \frac{\partial^2}{\partial z_iz_k} d(d-1)z_iz_k (1 + rz_1)^{d-2}\left|_{z=re_1}\right. = \inserted{(1+\delta_{ik})}d(d-1) (1+r^2)^{d-2},$
\item $\left\langle  \partial_{e_ie_k}, \partial_{e_je_l}\right\rangle_\poly = \frac{\partial^2}{\partial z_iz_k} d(d-1)z_jz_l (1 + rz_1)^{d-2}\left|_{z=re_1}\right. = 0$ when $\left\{i,k\right\} \ne \left\{j,l\right\}.$
\end{itemize}
\end{proof}

Let us now evaluate $\Psi$. Recall that
$$V_{re_1} = \left\{ f \in \poly \ : \ Df(re_1) = 0 \right\}$$
or, in other words, $f \in V_{re_1}$ if and only if
$$\innerP{f}{\partial_{e_i}}=0, \ 1 \le i \le n.$$
Thus, by lemma \ref{lem-2}-3, $\partial_{e_i}$, $\ 1 \le i \le n$, constitue an orthogonal basis of $V_{re_1}^\perp.$
We also have
$$\ker D^2 = \Span \left\{ \partial_{e_ie_j}, \ 1 \le i \le j \le n \right\}^\perp \cap V_{re_1}$$
hence,
\[
(\ker D^2)^\perp = \Span \left\{ P \partial_{e_ie_j}, \ 1 \le i \le j \le n \right\}
\]
where $P$ stands for the orthogonal projection onto $V_{r e_1}$.
We have seen that for $(i,j,k) \ne (1,1,1)$, $\partial{e_ie_j} \perp
\partial{e_k}$ (lemma \ref{lem-2}-5). Hence, 
\[
P \partial_{e_1e_1} \partial_{e_1e_1} - \partial_{e_1} 
\frac{ \innerP{\partial_{e_1e_1}}{\partial_{e_1}}}
{\| \partial e_1 \|_{\poly}^2 }
\]
and for $(i,j) \ne (1,1)$, 
$$P \partial_{e_ie_j} =  \partial_{e_ie_j}.$$
Let us now show that 
\[
\Psi (S) = \sum_{1 \le i \le j \le n}
S_{ij} \frac{ P \partial_{e_ie_j} }
{
\| P \partial_{e_ie_j} \|_{\poly}^2
}.
\]
Since this expression is clearly in $(\ker D^2)^\perp$ it suffices to prove that $D^2 \circ \Psi (S) = S$ for any $S \in \sym$ i.e.
$$D^2\Psi(S)(re_1)(e_k,e_l) = S_{kl}$$
or, using lemma \ref{lem-1}, that
$$\innerP{\partial_{e_ke_l}}{
\sum_{1 \le i \le j \le n}
S_{ij}
\frac{P\partial_{e_ie_j}}
{\left\| P\partial_{e_ie_j} \right\|_\poly^2}
} = S_{kl}.$$
This last equality holds because $P\partial_{e_ie_j}$, $1 \le i \le j \le n$, constitue an orthogonal basis of $(\ker D^2)^\perp$.

It is important to have in mind that $\Psi$ is not an isometry, we have
\[
\| \Psi (S)\|_\poly^2 = \sum_{1 \le i \le j \le n} \frac{S_{ij}^2}{\|P \partial_{e_ie_j}\|_{\poly}^2}
\]

We introduce now the functions
\[
A(d,r) = \sqrt{ \frac{ d(d-1)r^4 + 2dr^2 + 2}{(d r^2 + 1)R^4}}
\]
and
\[
B(d,r) = \sqrt{ \frac{ (d-1)r^2 + 1}{R^2}} \ ,
\]
where again $R = \sqrt{1+r^2}$.

\begin{lemma} \label{lem-3} Let $i \le j$. Then,
\[
\| P \partial_{e_i e_j} \|_{\poly} ^2
d(d-1) R^{2d-4} \times 
\left\{
\begin{array}{ll}
A(d,r)^2 & \text{if $i = 1$ and $j = 1$} \\
B(d,r)^2 & \text{if $i=1$ and $j \ne 1$} \\
{(1+\delta_{ij})} 
& \text{if $i \ne 1$ and $j \ne 1$}
\end{array}
\right.
\]
with $\delta_{ij}=1$ when $i=j$ and $0$ otherwise. 
\end{lemma}

Let us now compute $\det \Psi^*\Psi$. For any $f = \sum_{1 \le i \le j \le n} f_{ij} P \partial_{e_ie_j} \in (\ker D^2)^\perp$ and for any $S \in \sym$ we have
{\[
\left\langle \Psi^*(f),S\right\rangle = 
\left\langle f,\Psi(S)\right\rangle_\poly = 
\sum_{1 \le i \le j \le n} f_{ij} S_{ij}
\]
}
{Therefore, we have always for any $T \in \mathcal S_n$:
\[
\langle \Psi^* \Psi (T), S \rangle = 
\sum_{1 \le i \le j \le n} \frac{T_{ij}S_{ij}}
{\left\| P\partial_{e_ie_j} \right\|_\poly^2}
\]
We write
the matrix of the operator $\Psi^* \Psi$ with respect to the orthonormal
basis of $\mathcal S$ given by $\mathrm{e}_1 \mathrm{e}_1^T, \dots,
\mathrm{e}_n \mathrm{e}_n^T$ and then, for $i < j$, 
$\frac{1}{\sqrt{2}} \left(\mathrm{e}_i \mathrm{e}_j^T + 
\mathrm{e}_j \mathrm{e}_i^T)\right)$: 
\[
\Psi^* \Psi \left[
\begin{matrix} 
\frac{1}{\|P\partial_{e_1e_1}\|^2}  \\
& \ddots \\
& & \frac{1}{\|P\partial_{e_ne_n}\|^2}  \\
& & & \frac{1}{2 \|P\partial_{e_1e_2}\|^2}  \\
& & & & \ddots  \\
& & & & & \frac{1}{2 \|P\partial_{e_{n-1}e_n}\|^2} \\
\end{matrix} \right]
\]
}
Using lemma \ref{lem-2} we obtain:

\begin{lemma} \label{lem-4}
\[
(\det \Psi^* \Psi)^{1/2} = 
{2^{-\frac{(n+2)(n-1)}{4}}}
\left( d(d-1) R^{2d-4}\right)^{-\frac{n(n+1)}{4}} A(d,r)^{-1} B(d,r)^{-(n-1)}.
\]
\end{lemma}
At this point

\begin{proposition}\label{pro-5}
Under the conditions above,
{
\[
E_U
\frac{\alpha_n}{d^{n/2}}
\int_{0}^{\infty}
\frac{\left( \det \Psi^* \Psi \right)^{\frac{1}{2}} r^{n-1}\dd r}{
(dr^2+1)^{1/2}
R^{(d-1)n-1}
}
\int_{D^2(U \cap V_{re_1})}
\frac{ \left| \det S \right| }
{(2 \pi)^{\dim \sym / 2}} e^{-\|\Psi(S)\|_\poly^2/ 2} \dd \sym
.\]
In particular, 
\[ C_{d,n}
\frac{\alpha_n}{d^{n/2}}
\int_{0}^{\infty}
\frac{\left( \det \Psi^* \Psi \right)^{\frac{1}{2}} r^{n-1}\dd r}{
(dr^2+1)^{1/2}
R^{(d-1)n-1}
}
\int_{\sym}
\frac{ \left| \det S \right| }
{(2 \pi)^{\dim \sym / 2}} e^{-\|\Psi(S)\|_\poly^2/ 2} \dd \sym
,\]
and
\[ E_{d,n}
\frac{\alpha_n}{d^{n/2}}
\int_{0}^{\infty}
\frac{\left( \det \Psi^* \Psi \right)^{\frac{1}{2}} r^{n-1}\dd r}{
(dr^2+1)^{1/2}
R^{(d-1)n-1}
}
\int_{\sym^{++}}
\frac{ \det S }
{(2 \pi)^{\dim \sym / 2}} e^{-\|\Psi(S)\|_\poly^2/ 2} \dd \sym
\]
}
where $\sym^{++}$ denotes the set of positive definite matrices. When $n=1$,
$$
C_{d,1} = 2 E_{d,1} = 
\frac{ 2\sqrt{d-1}}{\pi}
\int_{0}^{\infty}
\frac{\sqrt{ d(d-1) r^4 + 2 d r^2 + 2}}{
(dr^2+1)(r^2+1)} 
\dd r .$$
\end{proposition}

\begin{proof} The three first formulas are obtained in combining proposition \ref{pro-4}, equation \ref{intout} and lemma \ref{lem-4}. For the case $n=1$ we obtain
$$E_{d,1} = \frac{2}{ d \sqrt{d-1} \sqrt{2 \pi}}
\int_{0}^{\infty}
\frac{\dd r}{
A 
(dr^2+1)^{1/2}
R^{2d-4} 
}
\ 
\int_{0}^{\infty}
\frac{ s }
{\sqrt{2 \pi}} e^{-\frac{ s^2}{2 d(d-1)R^{2d-4}A^2}} \dd s = $$
$$\frac{\sqrt{d-1}}{\pi}
\int_{0}^{\infty}
\frac{\sqrt{ d(d-1) r^4 + 2 d r^2 + 2}}{
(dr^2+1)(r^2+1)} 
\dd r .$$
The identity $C_{d,1} = 2 E_{d,1}$ is easy. \end{proof}

\section{Some integral lemmas}
{The term $e^{-\|\Psi(S)\|_{\mathcal P}^2/2}$ in the inner integrals
of Proposition~\ref{pro-5} can be simplified through additional
changes of coordinates. We reparametrize the spaces $\mathcal S_n$
and $\mathcal S_{n}^{++}$ though a stretching $S \mapsto 
T=\Delta^{-1} S \Delta^{-1}$.}

{The stretching coefficients are} $\Delta_i = \left(2d(d-1)R^{2d-4}\right)^{1/4}$ for $i \ge 2$, $\Delta_1 = B(d,r) \Delta_2$ and, $\Delta = \mbox{Diag}(\Delta_1, \Delta_2, \ldots , \Delta_n)$. We obtain
$$
\| \Psi(S) \|_{\poly}^2 = \frac{1}{d(d-1)R^{2d-4}}\left(
\frac{S_{11}^2}{A^2} + \sum_{j=2}^n \frac{S_{1j}^2}{B^2} + \sum_{1 < i \le j \le n} {\frac{1}{1+\delta_{ij}}} S_{ij}^2
\right)
$$
and
{
\[
\| \Delta^{-1}S\Delta^{-1} \|^2 = \frac{1}{d(d-1)R^{2d-4}}\left(
\frac{S_{11}^2}{2B^4} + \sum_{j=2}^n \frac{S_{1j}^2}{B^2} + \sum_{1 < i \le j \le n} \frac{1}{1+\delta_{ij}}S_{ij}^2
\right)
\]
}
so that 
{
\[
\| \Psi(S) \|_{\poly}^2 = \| \Delta^{-1}S\Delta^{-1} \|^2 + \left(\frac{1}{A^2} - \frac{1}{2B^4}\right)\frac{S_{11}^2}{d(d-1)R^{2d-4}}.
\]
}
Let us define $T = \Delta^{-1}S\Delta^{-1}$. We get 
{
\[
\| \Psi(S) \|_{\poly}^2 = \| T \|^2 + \left( \frac{2 B^4}{A^2} - 1 \right)T_{11}^2
\]
}
so that, via this change of variable,
\[
\int_{D^2(U)} \frac{\left| \det S \right|}{\sqrt{2 \pi}^{\dim \sym}} 
e^{-\| \Psi(S) \|^2 / 2} \dd S = \]
\[
\left( \prod_{i=1}^n \Delta_i \right)^{ n+3 }
\int_{\Delta^{-1} D^2(U\inserted{\cap V_r}) \Delta^{-1}} 
\frac{\left| \det T \right|}{\sqrt{2 \pi}^{\dim \sym}} 
e^{-\frac{1}{2}\left(\| T\|^2  + \left( \frac{{2} B(d,r)^4}{A(d,r)^2} - 1\right)T_{11}^2 \right)} \dd T.
\]
{
If $U \subset V$, we define the auxiliary quantity
\[
C_U (d,r,n) = 
\int_{\Delta^{-1} D^2(U \cap V_r) \Delta^{-1}} 
\frac{\left| \det T \right|}{\sqrt{2 \pi}^{\dim \sym}} 
e^{-\frac{1}{2}\left(\| T\|^2  + \left( \frac{2 B(d,r)^4}{A(d,r)^2} - 1\right)T_{11}^2 \right)} \dd T.
\]

There are two cases of interest corresponding to $U=V$ for the average of critical points and $U = V_+$ for the average number
of local minima. The corresponding functions are denoted $C_V (d,r,n)$ and $C_{V_+} (d,r,n)$. 
}
Using proposition \ref{pro-5} we get (the proof is easy and left to the reader)

\begin{proposition} \label{pro-6}
{
\[
E_U
\frac{2 \sqrt{2} (d-1)^{n/2}}{\Gamma(n/2)}
\int_{0}^{\infty}
\frac{((d-1)r^2+1)^2 }{R^2 \sqrt{d(d-1)r^4+2dr^2+2} 
}
\frac{r^{n-1}}{R^{n-1}}
C_U(d,r,n) \dd r
\]
}
Moreover
$$C_V(d,r,n) = \int_{\sym} 
\frac{\left| \det T \right|}{\sqrt{2 \pi}^{\dim \sym}} 
e^{-\frac{1}{2}\left(\| T\|^2  + \left( \frac{2 B(d,r)^4}{A(d,r)^2} - 1\right)T_{11}^2 \right)} \dd T$$
and
$$C_{V_+}(d,r,n) = \int_{\sym^{++}} 
\frac{ \det T }{\sqrt{2 \pi}^{\dim \sym}} 
e^{-\frac{1}{2}\left(\| T\|^2  + \left( \frac{2 B(d,r)^4}{A(d,r)^2} - 1\right)T_{11}^2 \right)} \dd T.$$
\end{proposition}

\section{Proof of Theorem \ref{general}}

To prove our main theorem we use both the proposition \ref{pro-6} and the case $d=2$ already investigated in the proposition \ref{pro-0}. We have
$$1 = C_{2,n} = \frac{2 \sqrt{2}}{\Gamma(n/2)} \int_0^\infty \frac{r^{n-1}}{R^{n-1}} \frac{C_V(2,r,n)}{\sqrt{2}} \dd r$$
and 
$$ C_V(2,r,n) = \int_{\sym} \frac{\left| \det T \right|}{(2 \pi)^{n(n+1)/4}} e^{-\frac{1}{2}\left(\left\|T\right\|^2 + 2r^2 T_{11}^2 \right)} \dd T.$$
%The rational expression appearing in $C_{d,n}$ is bounded by
%$$\frac{1}{\sqrt{2}} \le \frac{((d-1)r^2+1)^2 }{R^2 \sqrt{d(d-1)r^4+2dr^2+2}} \le \frac{(d-1)^{3/2}}{d^{1/2}}$$
%thus
%$$C_{d,n} \le \frac{2 \sqrt{2} (d-1)^{n/2}}{\Gamma(n/2)} \int_{0}^{\infty}\frac{(d-1)^{3/2}}{d^{1/2}}\frac{r^{n-1}}{R^{n-1}} C_V(d,r,n) \dd r.$$
%From the inequality (easy to obtain)
%$$2r^2 \le 2(d-1)r^2 \le \frac{2 B(d,r)^4}{A(d,r)^2} - 1 \le \frac{5}{4}2(d-1)r^2$$
%we get $C_V(d,r,n) \le C_V(2,r,n)$ so that
%$$C_{d,n} \le \frac{2 \sqrt{2} (d-1)^{n/2}}{\Gamma(n/2)} \int_{0}^{\infty}\frac{(d-1)^{3/2}}{d^{1/2}}\frac{r^{n-1}}{R^{n-1}} C_V(2,r,n) \dd r =$$
%$$(d-1)^{n/2} \frac{(d-1)^{3/2}}{d^{1/2}} \sqrt{2}
%$$
%and we are done. The case of local minima is treated via the same argument. 

\begin{lemma} \label{scaling} The quantity
$\Lambda(d,r) = \frac{2 B(d,r)^4}{A(d,r)^2}-1$ satisfies, for all
$r > 0$ and $d \ge 2$, the scaling law:
\[
\Lambda(2, r\sqrt{d-1} ) \le 
\Lambda(d, r) \le 
\Lambda(2, \frac{\sqrt{5}}{2} r\sqrt{d-1} )  
\]
\end{lemma}

\begin{proof}
We write
\[
\Lambda(d,r) = 2 (d-1) r^2 + \frac{d-2}{d} \frac{(d-1)r^4}
{(d-1)r^4 + 2 r^ 2 + \frac{2}{d}}
\ . 
\]
The lower bound is now obvious. The upper bound is obtained
as follows:
\begin{eqnarray*}
\Lambda(d,r) &=& 2 (d-1) r^2 + \frac{d-2}{d} \frac{(d-1)r^4}
{(d-1)r^4 + 2 r^ 2 + \frac{2}{d}}
\\
&\le&
\Lambda(d,r) = 2 (d-1) r^2 + \frac{d-2}{2d} (d-1)r^2
\\
&\le&
\frac{5}{4} \Lambda(2, r^2\sqrt{d-1})
\\ &=&
\Lambda(2, 
\frac{\sqrt{5}}{2} 
r^2 \sqrt{d-1}).
\end{eqnarray*}

\end{proof}
It follows from Lemma~\ref{scaling} that
\begin {eqnarray*}
C_V(d,r,n) &=&
\int_{S_n} 
\frac{|\det T |}{\sqrt{2 \pi}^{\dim S_n}}
e^{ -\frac{1}{2} \left( \|T\|^2 + \Lambda(d,r) T_{11}^2 \right)}
\\
&\le&
\int_{S_n} 
\frac{|\det T |}{\sqrt{2 \pi}^{\dim S_n}}
e^{ -\frac{1}{2} \left( \|T\|^2 + \Lambda(2,r\sqrt{d-1}) T_{11}^2 \right)}
\\
&=&
C_V(2, r\sqrt{d-1}, n)
\end{eqnarray*}
and similarly $C_{V+}(d,r,n)\le C_{V+} (2, r\sqrt{d-1}, n)$. Now we have:
\begin{eqnarray*} 
C_{d,n} &=&
\frac{2 \sqrt{2} (d-1)^{n/2}}{\Gamma({n}/{2})}
\int_0^\infty 
\frac{((d-1)r^2 + 1)^ 2}{R^2 \sqrt{d(d-1)r^4 + 2dr^2 + 2}}
\frac{r^{n-1}}{R^{n-1}} 
C_V(d,r,n) \ dr \\
&\le& 
\frac{2 \sqrt{2} (d-1)^{n/2}}{\Gamma({n}/{2})}
\int_0^\infty 
\frac{((d-1)r^2 + 1)^ 2}{R^2 \sqrt{d(d-1)r^4 + 2dr^2 + 2}}
\frac{r^{n-1}}{R^{n-1}} 
C_V(2,r\sqrt{d-1},n) \ dr
\ .
\end{eqnarray*}
We set $s = r\sqrt{d-1}$ and $S=\sqrt{d-1+s^2}$ to obtain:
\begin{eqnarray*} 
C_{d,n} &\le&
\frac{2 \sqrt{2} (d-1)^{n/2}}{\Gamma({n}/{2})}
\int_0^\infty 
\frac{({d-1})s^ 2}{{S^2} \sqrt{\frac{d}{d-1} s^4 + 2\frac{d}{d-1}s^2 + 2}}
\frac{s^{n-1}}{S^{n-1}} \
C_V(2,s,n) \ \frac{1}{\sqrt{d-1}} \ ds.
\end{eqnarray*}
Since 
$$\frac{({d-1})s^ 2}{{S^2} \sqrt{\frac{d}{d-1} s^4 + 2\frac{d}{d-1}s^2 + 2}} \le \sqrt{\frac{d-1}{d}}$$
and $C_{2,n} = 1$ we obtain
$$C_{d,n} \le \frac{(d-1)^{(n+2)/2}\sqrt{2}}{\sqrt{d}}$$
and the same argument holds for $E_{d,n}$.

\section{The Riemann surface}

We rewrite the case $n=1$ (proposition \ref{pro-5}) for convenience as:
\begin{equation}\label{one4}
E_{d,1} = \frac{(d-1)\sqrt{d}}{2 \pi} 
\int_{\mathbb R} 
g(z) \dd z
\end{equation}
with
\[
g(z) = \frac{ \sqrt{z^4+\frac{2}{d-1}z^2+\frac{2}{d(d-1)}} }{(1+z^2)(1 + dz^2)}
\]

At this point we encounter a classical situation: we want to compute a line
integral of a function g(z), which is a two-branched meromorphic function
of $\mathbb C$. In order to apply the residue theorem, we need first to
replace $g$ by a regular meromorphic function, defined in the relevant
Riemann surface $R$. The branching points of the Riemann surface are
the roots of the polynomial inside the square root. If we set

\[
\zeta = \sqrt{\frac{-1 + i \sqrt{1 - \frac{2}{d}}}{d-1}} 
\]
with the branch of the external square root in such a way that
$\zeta$ belongs to the positive quadrant, 
we can now factorize
\[
z^4+\frac{2}{d-1}z^2+\frac{2}{d(d-1)} (z-\zeta)(z-\bar \zeta)(z+\zeta)(z+\bar \zeta) 
\ . \]

It follows that the Riemann Surface $R$
is a twofold cover of $\mathbb C$ with
branch points $\zeta$, $-\bar \zeta$, $-\zeta$, $\bar \zeta$.

Let $\gamma$ be the arc of circle (centered in the origin)
joining $-\bar \zeta$ to $\zeta$ crossing
the positive imaginary axis. Notice that it croses
the segment $[i/\sqrt{d}, i]$. Let $\mathcal D$ denote the upper half plane
with $\gamma$ removed.

Then, the positive branch of $\sqrt{z^4+\frac{2}{d-1}z^2+\frac{2}{d(d-1)}}$
on $\mathbb R$ extends to a unique branch on $\mathcal D$. The square
root is real and positive on $[0,i|\zeta|]$ and real and negative on
$[i|\zeta|, i \infty)$. 

The residue theorem is now:

\[
\int_{\mathbb R} g(z) \dd z
- 2 \int_{\gamma} g(z) \dd z
2 \pi i \res_{[z=i/\sqrt{d}]} g(z) + 2 \pi i \res_{[z=i]} g(z) 
\]

Residues are respectively $\frac{-i}{2 (d-1) \sqrt{d}}$
and $\frac{-i \sqrt{d-2}}{2 (d-1) \sqrt{d}}$. Therefore,
\begin{equation}\label{one5}
E_{d,1} = \frac{1}{2} + \frac{\sqrt{d-2}}{2} +  
 \frac{(d-1)\sqrt{d}}{\pi} 
 \int_{\gamma} g(z) \dd z
\end{equation}

(We mean the integral of the branch that is positive on $i |\zeta|$).

Now, in order to integrate $g(z)$, we introduce a linear fractional
transformation mapping the real line onto the circle containing $\gamma$.
Namely,

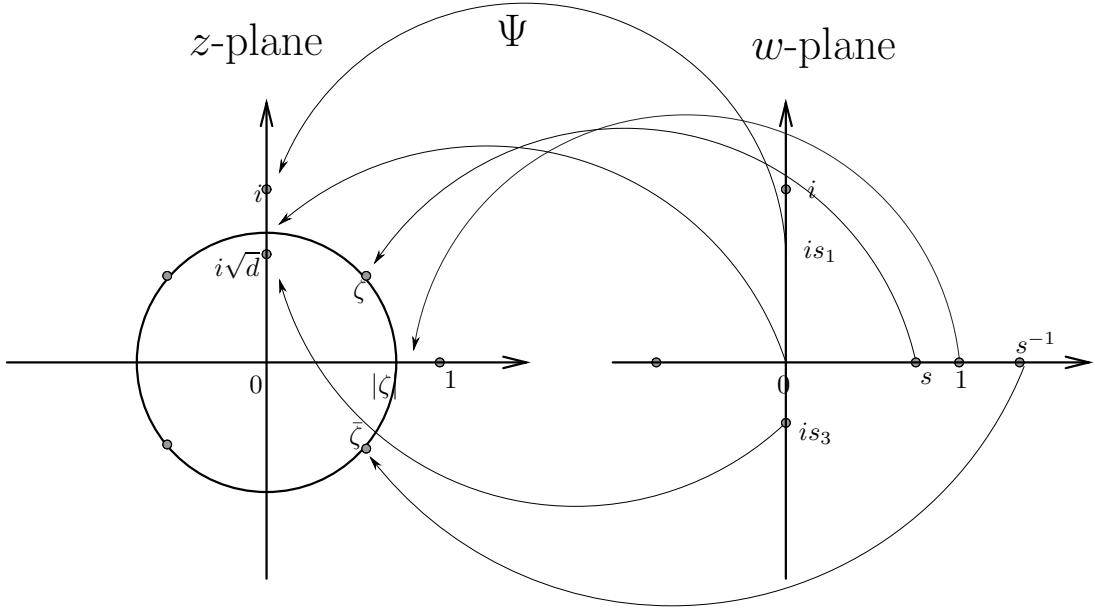
\begin{figure}
\setlength{\unitlength}{0.00050in}
\begingroup\makeatletter\ifx\SetFigFont\undefined
% extract first six characters in \fmtname
\def\x#1#2#3#4#5#6#7\relax{\def\x{#1#2#3#4#5#6}}%
\expandafter\x\fmtname xxxxxx\relax \def\y{splain}%
\ifx\x\y   % LaTeX or SliTeX?
\gdef\SetFigFont#1#2#3{%
  \ifnum #1<17\tiny\else \ifnum #1<20\small\else
  \ifnum #1<24\normalsize\else \ifnum #1<29\large\else
  \ifnum #1<34\Large\else \ifnum #1<41\LARGE\else
     \huge\fi\fi\fi\fi\fi\fi
  \csname #3\endcsname}%
\else
\gdef\SetFigFont#1#2#3{\begingroup
  \count@#1\relax \ifnum 25<\count@\count@25\fi
  \def\x{\endgroup\@setsize\SetFigFont{#2pt}}%
  \expandafter\x
    \csname \romannumeral\the\count@ pt\expandafter\endcsname
    \csname @\romannumeral\the\count@ pt\endcsname
  \csname #3\endcsname}%
\fi
\fi\endgroup
\centerline{
{%\renewcommand{\dashlinestretch}{30}
\begin{picture}(11294,6303)(0,-10)
\put(7081.018,2256.092){\arc{5720.609}{3.2896}{6.1669}}
\blacken\path(4243.090,2801.372)(4252.000,2678.000)(4302.196,2791.050)(4266.450,2760.748)(4243.090,2801.372)
\put(5012.897,1520.420){\arc{6545.897}{3.9933}{5.9654}}
\blacken\path(2929.263,4083.390)(2857.000,3983.000)(2968.004,4037.573)(2921.143,4037.237)(2929.263,4083.390)
\put(6440.730,1871.536){\arc{6209.497}{3.7235}{6.0652}}
\blacken\path(3890.144,3693.925)(3847.000,3578.000)(3939.621,3659.983)(3894.518,3647.268)(3890.144,3693.925)
\put(6931.682,3951.243){\arc{7885.805}{0.3775}{2.4877}}
\blacken\path(3900.129,1477.696)(3802.000,1553.000)(3853.147,1440.377)(3854.247,1487.226)(3900.129,1477.696)
\put(5401.987,3556.607){\arc{5444.586}{3.5045}{6.2423}}
\blacken\path(2874.438,4645.458)(2857.000,4523.000)(2930.013,4622.845)(2888.658,4600.806)(2874.438,4645.458)
\put(5933.825,4230.834){\arc{6375.097}{0.8142}{2.8772}}
\blacken\path(2919.263,3291.120)(2857.000,3398.000)(2861.641,3274.394)(2880.416,3317.330)(2919.263,3291.120)
\thicklines
\put(2722,2543){\ellipse{2700}{2700}}
\thinlines
\texture{55888888 88555555 5522a222 a2555555 55888888 88555555 552a2a2a 2a555555 
	55888888 88555555 55a222a2 22555555 55888888 88555555 552a2a2a 2a555555 
	55888888 88555555 5522a222 a2555555 55888888 88555555 552a2a2a 2a555555 
	55888888 88555555 55a222a2 22555555 55888888 88555555 552a2a2a 2a555555 }
\put(4522,2543){\shade\ellipse{90}{90}}
\put(4522,2543){\ellipse{90}{90}}
\put(2722,4343){\shade\ellipse{90}{90}}
\put(2722,4343){\ellipse{90}{90}}
\put(8122,4343){\shade\ellipse{90}{90}}
\put(8122,4343){\ellipse{90}{90}}
\put(9922,2543){\shade\ellipse{90}{90}}
\put(9922,2543){\ellipse{90}{90}}
\put(9472,2543){\shade\ellipse{90}{90}}
\put(9472,2543){\ellipse{90}{90}}
\put(10552,2543){\shade\ellipse{90}{90}}
\put(10552,2543){\ellipse{90}{90}}
\put(3757,3443){\shade\ellipse{90}{90}}
\put(3757,3443){\ellipse{90}{90}}
\put(1687,3443){\shade\ellipse{90}{90}}
\put(1687,3443){\ellipse{90}{90}}
\put(3757,1643){\shade\ellipse{90}{90}}
\put(3757,1643){\ellipse{90}{90}}
\put(1687,1688){\shade\ellipse{90}{90}}
\put(1687,1688){\ellipse{90}{90}}
\put(6772,2543){\shade\ellipse{90}{90}}
\put(6772,2543){\ellipse{90}{90}}
\put(2722,3668){\shade\ellipse{90}{90}}
\put(2722,3668){\ellipse{90}{90}}
\put(8122,1913){\shade\ellipse{90}{90}}
\put(8122,1913){\ellipse{90}{90}}
\thicklines
\path(22,2543)(5422,2543)
\path(5182.000,2483.000)(5422.000,2543.000)(5182.000,2603.000)
\path(2722,293)(2722,5243)
\path(2782.000,5003.000)(2722.000,5243.000)(2662.000,5003.000)
\path(8122,293)(8122,5243)
\path(8182.000,5003.000)(8122.000,5243.000)(8062.000,5003.000)
\path(6322,2543)(11272,2543)
\path(11032.000,2483.000)(11272.000,2543.000)(11032.000,2603.000)
\put(9877,2273){\makebox(0,0)[lb]{\smash{{{\SetFigFont{10}{20.4}{rm}$1$}}}}}
\put(9517,2318){\makebox(0,0)[lb]{\smash{{{\SetFigFont{10}{20.4}{rm}$s$}}}}}
\put(10507,2633){\makebox(0,0)[lb]{\smash{{{\SetFigFont{10}{20.4}{rm}$s^{-1}$}}}}}
\put(8032,2228){\makebox(0,0)[lb]{\smash{{{\SetFigFont{10}{20.4}{rm}$0$}}}}}
\put(8257,1733){\makebox(0,0)[lb]{\smash{{{\SetFigFont{10}{20.4}{rm}$is_3$}}}}}
\put(8302,3623){\makebox(0,0)[lb]{\smash{{{\SetFigFont{10}{20.4}{rm}$is_1$}}}}}
\put(8347,4253){\makebox(0,0)[lb]{\smash{{{\SetFigFont{10}{20.4}{rm}$i$}}}}}
\put(2587,4208){\makebox(0,0)[lb]{\smash{{{\SetFigFont{10}{20.4}{rm}$i$}}}}}
\put(4567,2273){\makebox(0,0)[lb]{\smash{{{\SetFigFont{10}{20.4}{rm}$1$}}}}}
\put(3622,3218){\makebox(0,0)[lb]{\smash{{{\SetFigFont{10}{20.4}{rm}$\zeta$}}}}}
\put(3577,1688){\makebox(0,0)[lb]{\smash{{{\SetFigFont{10}{20.4}{rm}$\bar \zeta$}}}}}
\put(2542,2228){\makebox(0,0)[lb]{\smash{{{\SetFigFont{10}{20.4}{rm}$0$}}}}}
\put(2182,3443){\makebox(0,0)[lb]{\smash{{{\SetFigFont{10}{20.4}{rm}$i\sqrt{d}$}}}}}
\put(3802,2228){\makebox(0,0)[lb]{\smash{{{\SetFigFont{10}{20.4}{rm}$|\zeta|$}}}}}
\put(5107,5828){\makebox(0,0)[lb]{\smash{{{\SetFigFont{17}{20.4}{rm}$\Psi$}}}}}
\put(7762,5693){\makebox(0,0)[lb]{\smash{{{\SetFigFont{17}{20.4}{rm}$w$-plane}}}}}
\put(1912,5738){\makebox(0,0)[lb]{\smash{{{\SetFigFont{17}{20.4}{rm}$z$-plane}}}}}
\end{picture}
}}
\caption{The Linear Fractional Map $w \mapsto z=\Psi(w)$}.
\end{figure}

\[
\Psi(w) = 
\frac{Aw + B}{Cw + D} 
\]

with $A = |\zeta|$, $ B = i |\zeta|$, $C=i$, $D=1$.
For the record, $AD-BC = 2 |\zeta|$

Let $s = \frac{\real(\zeta)}{|\zeta| + \imag(\zeta)}$. 
Define also $s_1 = \frac{1-|\zeta|}{1+|\zeta|}$, $s_2=s_1^{-1}$,
$s_3= \frac{1-|\zeta|\sqrt{d}}{1+|\zeta|\sqrt{d}}$ and $s_4=s_3^{-1}$.
We have
the following mapping table for $\Psi$:
\medskip
\par

\centerline{
\begin{tabular}{||c|c||c|c||}
\hline
\hline
$w$ & $\Psi(w)$ 
& $w$ & $\Psi(w)$ 
\\
\hline
\hline
$-1$ & $-|\zeta|$      & $i s_1$       & $i$   \\
$0$  & $i |\zeta|$     & $i s_2$      & $-i$  \\
$1$  & $|\zeta|$       & $i s_3$& $i/\sqrt{d}$ \\
$-s^{-1}$ & $-\zeta$      &  $i s_4$ & $-i/\sqrt{d}$\\
$-s$ & $- \bar \zeta$    &&   \\
$s$  & $\zeta$           &&   \\
$s^{-1}$ & $\bar \zeta$  &&   \\
\hline
\hline
\end{tabular}}

\medskip
\par

Changing coordinates,
\[
\int_{\gamma}
g(z) \dd z
= 2 c(d)
\real \int_{[0,s]}
\frac{\sqrt{ (w^2 - s^2) (w^2 - s^{-2}) }}
{
\prod_{k=1}^4 (w-i s_k)
}
\
\dd w
\]
with
\[
c(d)= \frac{ (AD-BC) \sqrt{A^4+ \frac{2}{d-1}A^2 C^2 + \frac{2}{d(d-1)}C^4}}
{(A^2 + C^2) (dA^2 + C^2)}
\in
O(d^{-3/2}) 
\]

(More precisely: $\lim d^{3/2} C(d) = -2^{7/4} \frac{\sqrt{2-\sqrt{2}}}{\sqrt{2}-1}
\simeq -6.2151$).

At this point, good practice seems to be:
\begin{enumerate}
\item Multiply numerator and denominator by the conjugate of the denominator,
in order to obtain a real polynomial in the denominator.
\item Multiply numerator and denominator by the square root.
\item Expand in partial fractions.
\item Put into Legendre normal form.
\item Write down the integral in terms of elliptic functions $K$ and $\Pi$.
\end{enumerate}

We expand the integrand in partial fractions:

\begin{eqnarray*}
\int_{\gamma}
g(z) \dd z
&=& 
2 c(d)
\int_{[0,s]}
\frac{1}{\sqrt{ (w^2 - s^2) (w^2 - s^{-2}) }}
\left(
1
+
\sum_{k=1}^4 \real \frac{R_k}{w-i s_k}
\right)
\
\dd w
\\
&=& 
2 c(d)
\int_{[0,s]}
\frac{1}{\sqrt{ (w^2 - s^2) (w^2 - s^{-2}) }}
\left(
1
+
\sum_{k=1}^4 \frac{s_k^{-2} \real \left(R_k (w + i s_k)\right)}{1+w^2 s_k^{-2}}
\right)
\
\dd w
\\
&=& 
2 c(d)
\int_{[0,s]}
\frac{1}{\sqrt{ (w^2 - s^2) (w^2 - s^{-2}) }}
\left(
1
+
\sum_{k=1}^4 \frac{s_k^{-1} R_k i}{1+w^2 s_k^{-2}}
\right)
\
\dd w
\end{eqnarray*}
(the last step uses the fact that all residues $R_k$ are pure imaginary).
Residues are given in Table~\ref{arguments}.
We use formula \cite{abr} [17.4.45] to compute the {\em parameter} $m=s^4$.
Then we set $\sin \alpha = s^{2}$ above, and also $w = s \sin \theta$ to
obtain the Legendre normal form:

\[
\int_{\gamma}
g(z) \dd z
= 2 c(d) s
\int_{0}^{\pi/2}
\frac{1}{\sqrt{ 1 - \sin^2 \alpha \sin^2 \theta}}
\left(
1
+
\sum_{k=1}^4 \frac{R_k s_k i} {1-n_k \sin^2 \theta}
\right)
\
\dd \theta
\]
This is a combination of one complete elliptic integral of the first kind
and 4 complete elliptic integrals of the third kind.
The {\em arguments} $n_k=-s^2 s_k^{-2}$ of the integrals of the 
third kind are given 
in Table~\ref{arguments}
\medskip
\par
\begin{table}
\centerline{\begin{tabular}{||l|c|c||}
\hline
\hline
Pole & Residue $R_k$ & Argument $n_k$ \\
\hline
$i s_1$ & 
$ i\frac{(s_1-s^2)(s_1-s^{-2})} { (s_1-s_2)(s_1-s_3)(s_1-s_4)} $
& $-\frac{s^2}{s_1^2}$    \\
$i s_2$ &     
$ i\frac{(s_2-s^2)(s_2-s^{-2})} { (s_2-s_1)(s_2-s_3)(s_2-s_4)} $
& $-\frac{s^2}{s_2^2}$    \\
$i s_3$ &     
$ i\frac{(s_3-s^2)(s_3-s^{-2})} { (s_3-s_1)(s_3-s_2)(s_3-s_4)} $
& $-\frac{s^2}{s_3^2}$    \\
$i s_4$ &     
$ i\frac{(s_4-s^2)(s_4-s^{-2})} { (s_4-s_1)(s_4-s_2)(s_4-s_3)} $
& $-\frac{s^2}{s_4^2}$    \\
\hline
\hline
\end{tabular}}
\caption{\label{arguments}Residues and arguments}
\end{table}
\medskip
\par

Therefore,

\[
\int_{\gamma}
g(z) \dd z
= 2 c(d)
\left(
K(m) + 
\sum_{k=1}^4 R_k s_k i \Pi(n_k; m)
\right)
\]
where $K$ and $\Pi$ denote the complete elliptic integrals of the first
and third kind, respectively.

ASYMPTOTICS: $s \rightarrow \sqrt{2}-1$, so $m \rightarrow 0.029437251$,
$\alpha \rightarrow 0.172425997$ rad $\simeq 9º52'45.42''$.
Also, $s_1, s_2 \rightarrow 1$ and $s_3 = s_4^{-1} = \rightarrow (1-\sqrt{2})/(1+\sqrt{2})$.

EXPERIMENTAL DATA: The hypergeometric functions were evaluated using Romberg iteration. Coefficients and residues obtained symbolically and then numerically. Digits are not guaranteed to be all significative. 

\medskip 
\par
\centerline{\begin{tabular}{||r|l||}
\hline
\hline
$d$ & $E_U - \frac{1}{2} - \frac{\sqrt{d-2}}{2}$ \\
\hline
$3$ & $- 0.280134$ \\
$4$ & $- 0.319279$ \\
$5$ & $- 0.337448$ \\
$6$ & $- 0.348064$\\
$7$ & $- 0.355053$\\
$8$ & $- 0.360010$\\
$9$ & $- 0.363712$\\
$10$& $- 0.366583$ \\
$10^2$ & $- 0.387335$ \\
$10^3$ & $- 0.389199$ \\
$10^4$ & $- 0.389384$ \\
$10^5$ & $- 0.389402$ \\
$10^6$ & $- 0.389404$\\
$10^7$ & $- 0.389405$\\
$10^8$ & $- 0.389405$ \\
$10^9$ & $- 0.389405$\\
\hline
\hline
\end{tabular}}
\medskip 
\par

\begin{remark}
Mark Rybowicz \cite{ryb} provided the following alternative formula
for $C_{d,1} = 2 E_{d,1}$:
\begin{eqnarray*}
C_{d,1} =
&-& \frac{4d(u-2)}{\sqrt{u}(u-1)(u-d)\pi} K(v)
\\ &+& \frac{u+1}{\sqrt{u} (u-1)\pi } \Pi(-\frac{(u-1)^2}{4u}, v)
\\ &+& \frac{(2-d)(u+d)}{\sqrt{u} (d-u) \pi} \Pi(-\frac{ (d-u)^2} {4du}, v)
\end{eqnarray*}
where
\[
u=\sqrt{ \frac{2d}{d-1} }
\text{   and   }
v=\frac{\sqrt{2-u}}{2}
\]
His formula agrees with ours up to six decimal places.
\end{remark}

\section{Appendix: Asymptotics for $P_n$}

Let us recall briefly the large deviation principle for large radom matrices. A good reference is Guionnet \cite{gui} where we have taken most of the following. 

Let $X$ be a real $n \times n$ symmetric matrix. Its entries are independent Gaussian random variables with mean $0$ and variances $2/n$ for a diagonal entry and $1/n$ for an off-diagonal one. Thus $\sqrt{\frac{n}{2}} X$ is in $GOE(n)$. 

Let ${\cal P}(\R)$ denotes the set of probability measures on $\R$ equipped with the weak topology. For any $\lambda \in \R^n$ let 
$\delta(\lambda) = \frac{1}{n} \sum_{i=1}^n\delta_{\lambda_i}$ be the probability measure on $\R$ defined by
$$\delta(\lambda)(A) = \frac{1}{n} \# \left\{ i \ : \ \lambda_i \in A \right\}.$$
We also need the function $I : {\cal P}(\R) \rightarrow [0, \infty]$ defined by
$$I(\mu) = \frac{1}{4} \int x^2 d\mu(x) - \frac{1}{2} \int \int \log(\left| x-y \right|) d\mu(x)d\mu(y) - \frac{3}{8}$$
and the probability measure $Q^n$ on $\R^n$ with density
$$Q^n =  \frac{1}{Z^n} \prod_{1 \le i < j \le n}\left| \lambda_i - \lambda_j \right| \exp\left( - \frac{n}{2} \sum_{i=1}^n \lambda_i^2\right)d\lambda$$
with respect to the Lebesgue measure on $\R^n$ ($Z^n$ is the normalisation constant). Then, according to Guionnet \cite{gui} Theorem 3.1, for every closed set ${\cal A} \subset {\cal P}(\R)$ one has
$$\limsup_{n \rightarrow \infty} \frac{1}{n^2} \log Q^n\left(\left\{ \lambda \in \R^n : \delta(\lambda) \in {\cal A}\right\}\right) \le 
- \inf_{\mu \in {\cal A}} I(\mu).$$
To obtain an estimate on $P_n$ we take 
$${\cal A} = \left\{\mu \in  {\cal P}(\R) \ : \ \mu([0, \infty[) = 1\right\}$$
so that, when $\lambda$ is the vector of eigenvalues of the matrix $X$, we have $\mu \in {\cal A}$ if and only if $X$ is semi-positive definite or (almost surely) if and only if $X$ is positive definite. Since $Q^n$ is the joint law of the eigenvalues of $X$ we get
$$\limsup_{n \rightarrow \infty} \frac{1}{n^2} \log \mbox{Prob}\left\{ X \in \sym \ : \ X \mbox{ is pd} \right\} \le 
- \inf_{\mu \in {\cal A}} I(\mu)$$
where the probability is taken in $\sqrt{\frac{2}{n}}GOE(n)$. Since the set of positive definite matrices is invariant under scaling it is not too difficult to see that 
$$P_n = \mbox{Prob}\left\{ X \in \sym \ : \ X \mbox{ is pd} \right\}$$
so that
$$\limsup_{n \rightarrow \infty} \frac{1}{n^2} \log (P_n) \le 
- \inf_{\mu \in {\cal A}} I(\mu) = - \alpha.$$

It remains to explain why $\alpha$ is a positive number. The map $I$ is $\ge 0$, each sub-level set $\left\{ \mu : I(\mu) \le M \right\}$ is compact and, it achieves its minimum value at a unique probability measure on $\R$ described as Wigner's semicircular law 
$$\frac{1}{2\pi}\sqrt{4-x^2} dx$$
with support at $[-2,2]$. Since this measure is clearly not in our set ${\cal A}$ we get $\alpha > 0$ and we are done.

\end{document}